\documentclass[reqno,11pt]{amsart}
\usepackage{amsmath, latexsym, amsfonts, amssymb, amsthm, amscd}

\numberwithin{equation}{section}

\newtheorem{theorem}{Theorem}[section]

\newtheorem{proposition}[theorem]{Proposition}

\newtheorem{rem}[theorem]{Remark}

\newtheorem*{assumption}{Basic Hypothesis}

\newcommand{\bbE}{{\ensuremath{\mathbb E}} }

\newcommand{\bbP}{{\ensuremath{\mathbb P}} }

\newcommand{\bbZ}{{\ensuremath{\mathbb Z}} }

\newcommand{\cH}{{\ensuremath{\mathcal H}} }

\newcommand{\gb}{\beta}
\newcommand{\gG}{\Gamma}
\newcommand{\gd}{\delta}

\newcommand{\gep}{\varepsilon}       

\newcommand{\gl}{\lambda}
\newcommand{\gL}{\Lambda}

\newcommand{\gs}{\sigma}

\newcommand{\gp}{\varphi}

\newcommand{\go}{\omega}

\renewcommand{\tilde}{\widetilde}          
\DeclareMathSymbol{\leqslant}{\mathalpha}{AMSa}{"36} 
\DeclareMathSymbol{\geqslant}{\mathalpha}{AMSa}{"3E} 
\DeclareMathSymbol{\eset}{\mathalpha}{AMSb}{"3F}     
\newcommand{\dd}{\text{\rm d}}             

\newcommand{\R}{\mathbb{R}}
\newcommand{\Z}{\mathbb{Z}}
\newcommand{\N}{\mathbb{N}}

\DeclareMathOperator{\sign}{sign}

\def\bP{\ensuremath{\bs{\mathrm{P}}}} 
\def\bE{\ensuremath{\bs{\mathrm{E}}}}

\newcommand{\ind}{\bs{1}}

\def\bs{\boldsymbol}

\def\free{{f}}
\def\freea{\tilde {f}}

\title[Annealed bounds for pinning models]{
On  constrained annealed bounds \\
for pinning and wetting models}

\author{Francesco Caravenna}

\address{Universit\`a di Milano-Bicocca, Dipartimento di Matematica e Applicazioni, \mbox{Edificio} U5,\break
via Cozzi 53, 20125 Milano, Italy
\hfill \break
\indent \textit{and} \hfill\break
\indent Laboratoire de Probabilit{\'e}s de P 6 \& 7 (CNRS U.M.R. 7599)
and  Universit{\'e} Paris 7 -- Denis Diderot,
U.F.R. Mathematiques, Case 7012, 2 place Jussieu, 75251 Paris cedex 05, France}

\email{f.caravenna\@@sns.it}

\urladdr{http://www.matapp.unimib.it/\raisebox{0.11ex}{\tiny$\sim$}fcaraven/}

\author{Giambattista Giacomin}

\address{Laboratoire de Probabilit{\'e}s de P 6\ \& 7 (CNRS U.M.R. 7599) and  Universit{\'e} Paris 7 -- Denis Diderot, U.F.R. Mathematiques, Case 7012, 2 place Jussieu, 75251 Paris cedex 05, France}

\email{giacomin\@@math.jussieu.fr}

\urladdr{http://www.proba.jussieu.fr/pageperso/giacomin/GBpage.html}

\date{September 8, 2005}



\begin{document}

\begin{abstract}
The  free energy of quenched disordered systems
is bounded above by the free energy of the corresponding
annealed system. This bound may be improved by
applying the annealing procedure, which is just Jensen inequality,
after having modified the Hamiltonian in a way that
the quenched expressions are left unchanged.
This procedure is often viewed as a partial annealing
or as a constrained annealing, in the sense that
the term that is added may be interpreted as a
Lagrange multiplier on the disorder variables.

 In this note
we point out that, for a family of models, some of
which have attracted much attention,
the multipliers of the form of empirical averages of local functions
cannot improve on the basic annealed bound from the viewpoint
of characterizing the phase diagram.
This class of multipliers is the one that is suitable
for computations and it is often believed that in this
class one can  approximate arbitrarily well the quenched free energy.
\end{abstract}

\maketitle

\section{The framework and the main result}
\label{sec:intro}

\subsection{The set--up (I): linear chain models}
\label{sec:subintro} A number of disordered models of linear chains
undergoing localization or pinning effects can be put into the
following general framework. Let  $S:=\left\{ S_n
\right\}_{n=0,1,\ldots}$ be a process with $S_n$ taking values in
$\Z^d$, $d \in \N :=\{1,2,\ldots\}$ and law $\bP$.

The {\sl disorder} in the system is given by a sequence $\go
:=\left\{\go _n\right\}_n$ of IID random variables of law $\bbP$,
with~$\go_n$ taking values in~$\Gamma \subseteq \R$. As a matter of
fact we could simply set $\Gamma=\R$, however several examples that
we will present deal with the case in which $\Gamma$ is a finite set
and in this situation our results require no measurability
conditions. The disorder acts on the paths of $S$ via an Hamiltonian
that, for a system of size $N$, is a function $H_{N,\go}$ of the
trajectory $S$, but depending only on $S_0, S_1, \ldots, S_N$. One
is interested in the properties of the probability measures $\bP_{N,
\go}$ defined by giving the density with respect to $\bP$:
\begin{equation}
\label{eq:1} \frac{\dd \bP_{N, \go}}{\dd\bP} \left( S \right) \, =\,
\frac1{Z_{N, \go}} \exp\left(H_{N, \go}\left(S\right)\right),
\end{equation}
where $Z_{N, \go}:= \bE \left[\exp\left(H_{N,
\go}\left(S\right)\right)\right]$ is the normalization constant. Our
attention focuses on the asymptotic behavior   of $\log Z_{N, \go}$.

\smallskip

In the sequel we will assume:

\medskip
\begin{assumption} \rm
There exists a sequence $\left\{ D_n \right\}_n$ of subsets of
$\Z^d$ such that $\bP ( S_n \in D_n \text{ for } n=1,2, \ldots,
N)\stackrel{N\to \infty}{\asymp} 1 $, namely
\begin{equation}
\lim_{N\to \infty} \frac 1N \log \bP \left( S_n \in D_n \text{ for }
n=1,2, \ldots, N\right)=0,
\end{equation}
and such that $H_{N, \go}(S)=0$ if $S_n \in D_n$ for $n=1, 2,
\ldots, N$.
\end{assumption}
\medskip

One sees directly that this hypothesis implies
\begin{equation}
\label{eq:first} \liminf_{N \to \infty} \frac 1N \log Z_{N, \go} \,
\ge \lim_{N \to \infty} \frac 1N \log \bP \left( S_n \in D_n \text{
for } n=1,2, \ldots, N\right)\, =\, 0,
\end{equation}
$\bbP (\dd \go)$--a.s.. We will assume that $\left\{(1/N)\log Z_{N,
\go}\right\}_N$ is a sequence of integrable random variables that
converges
 in the $L^1\left( \bbP (\dd \go)\right)$ sense and
$\bbP(\dd \go)$--almost surely to a constant, the {\sl free energy},
that we will call $\free$. These assumptions are verified in the
large majority of the interesting situations, for example whenever
super/sub--additivity tools are applicable.

Of course \eqref{eq:first} says that $\free \ge 0$ and one is lead
to the natural question of whether $\free =0$ or $\free >0$. In the
instances that we are going to consider the free energy may be zero
or positive according to some parameters from which $H_{N, \go}(S)$
depends: $\free =0$ and $\free >0$ are associated to sharply
different behaviors of the system.

\medskip

In order to establish upper bounds on $\free$ one may apply directly
Jensen inequality ({\sl annealed bound}) obtaining
\begin{equation}
\label{eq:J}
\begin{split}
\free \;& =\; \lim_{N\to\infty} \, \frac 1N \; \bbE \big[ \log Z_{N,\go} \big] \\
& \le \; \liminf_{N \to \infty} \, \frac 1N \, \log \bbE \big[
Z_{N,\go} \big]\; =: \; \freea \in [0, \infty],
\end{split}
\end{equation}
and,  in our context,
  if $\freea=0$ then $\free =0$.
The annealed bound
 may be improved by adding to
$H_{N, \go}(S)$ an integrable function $A_N: \Gamma^\N \to \R$ such
that $\bbE\left[ A_N (\go) \right]=0$: in fact $\free$ as defined in
the first line of~\eqref{eq:J} is unchanged by such transformation,
while the second line of~\eqref{eq:J} may depend on the choice of
$\{A_N\}_N$. We stress that not only $\free$ is left unchanged by
$H_{N, \go}(S) \to H_{N, \go}(S)+A_N (\go)$, but $\bP_{N,\go}$
itself is left unchanged (for every $N$). Notice moreover that the
{\sl optimal choice} $A_N (\go) = -\log Z_{N,\go} + \bbE\left[  \log
Z_{N,\go}\right]$ yields the equality in \eqref{eq:J}.

In the sequel when we refer to $\freea$ we mean that  $Z_{N,\go}$ is
defined with respect to $H_{N, \go}$ satisfying the Basic Hypothesis
(no $A_N$ term added).

\subsection{The result}
What we prove in this note is that
\medskip

\begin{proposition}
\label{th:main} If $\freea >0$ then for every local bounded
measurable function $F: \Gamma ^\N \longrightarrow \R$ such that
$\bbE \left[ F (\go) \right]=0$ one has
\begin{equation}
\label{eq:main} \liminf_{N \to \infty} \frac 1 N \log \bbE \bE
\left[ \exp \left( H_{N, \go}(S) + \sum_{n=0}^N F(\theta_n
\go)\right)\right]\, > \, 0,
\end{equation}
where $(\theta_n \go)_m= \go_{n+m}$.
\end{proposition}

\medskip

We can sum up this result by saying that when  $ \free =0$ but $
\freea >0$ it is of no use  modifying  the Hamiltonian by  adding
the empirical average of a (centered) local (bounded measurable)
function.

Notice that requiring $F(\cdot)$ to be bounded and measurable is
superfluous if $\gG$ is a finite set. From now on the reader should
read {\sl local} as a short--cut for {\sl local, measurable and
bounded}. We take this occasion also to observe that in principle
one should be able to extend the result in the direction  of
unbounded $F(\cdot)$ or  of non IID disorder: this however requires
additional assumptions and leads far from the spirit of this note.


\smallskip

On a mathematical level it is not obvious that the free energy may
be approximated via empirical averages of a local function of the
disorder, because we are playing with an exchange of limits (recall
the optimal choice of~$A_N$ above). But we remark that in the
physical literature the approach of approximating the free energy
via what can be viewed as a constrained annealed computation, the
term $\sum_{n=0}^N F(\theta_n \go)$ being interpreted as a Lagrange
multiplier, is often considered as an effective way of approximating
the quenched free energy. Here we mention in particular
\cite{cf:Morita} and \cite{cf:Kuhn} in which this point of view is
taken up in a systematic way: the aim is to approach the quenched
free energy by constrained annealing via local functions $F$ that
are more and more complex, the most natural example being linear
combinations of correlations of higher and higher order.
\smallskip

The proof of Proposition \ref{th:main} is based on the simple
observation that whenever $A_N$ is centered
\begin{multline}
\label{eq:argproof}
 \frac 1 N \log
\bbE \bE\left[ \exp \left( H_{N, \go}(S) + A_N(\go) \right)\right]\,
\ge
\\
 \frac 1 N \log
\bbE \left[ \exp \left( A_N(\go) \right)\right] + \frac 1 N \log \bP
\left( S_n \in D_n \text{ for } n=1, 2, \ldots , N \right)\, =:\,
Q_N+P_N.
\end{multline}
By hypothesis $P_N=o(1)$ so one has to consider the asymptotic
behavior of $Q_N$. If $\liminf_N Q_N >0$ there is nothing to prove.
So let us assume that $\liminf_{N}Q_N=0$: in this case the inferior
limit of the left--hand side of \eqref{eq:argproof} may be  zero and
we want to exclude this possibility when $\tilde f >0$ and $A_N
(\go)= \sum_{n=0}^N F(\theta_n \go)$, $F$ local and  centered (of
course in this case $\lim_N Q_N$ does exist). And in
Proposition~\ref{th:main2} below  in fact we show that if $\log
\bbE\left[ \exp\left(A_N(\go)\right)\right] =o(N)$, then $\sup_\go
\left\vert A_N (\go) \right \vert = o(N)$ and therefore
 the corresponding constrained annealing is just the standard annealing.

\medskip

\begin{rem}
\label{th:rem} \rm We stress that our Basic Hypothesis is more
general than it may look at first. As already observed, one has  the
freedom of adding to the Hamiltonian $H_{N, \go}(S)$ any term  that
does not depend on $S$ (but  possibly does depend on $\go$ and $N$)
without changing the model $ \bP_{N,\go}$.
 It may therefore happen
 that the {\sl natural} formulation of the Hamiltonian does not satisfy
 our Basic Hypothesis, but it does after a suitable additive correction.
  This  happens for example  in \S \ref{sec:cop} below:
  the additive correction in that case is linear in $\go$ and it corresponds to
  what in \cite{cf:ORW} is called {\sl first order} Morita approximation.
  In these terms, Proposition \ref{th:main} is saying that
  {\sl higher order} Morita approximations cannot improve
  the bound on the critical curve found with the first order computation.
\end{rem}

\medskip

\begin{rem}
\label{th:rem_added} \rm In the Morita approach of
\cite{cf:Kuhn,cf:Morita}, when applied to spin systems, it was also
taken for granted in the physics literature that the infinite volume
measure describing the joint distribution of disorder variables and
spin variables can be described  as  Gibbs measure with  a proper
(absolutely summable) Hamiltonian. This was shown to be false in
general, and potentials with weaker summability properties are
needed \cite{cf:add1,cf:add2}. This phenomenon underlines from a
different perspective that local dependence of the {\sl Morita
potential} on the disorder variables {\sl is not enough}.
\end{rem}

\medskip

Let us now look at applications of Proposition \ref{th:main}.

\subsubsection{Random rewards or penalties at the origin}
Let $S$, $S_0=0\in \Z ^d$, be a random walk with centered IID non
degenerate increments $\{X_n\}_n$, $(X_n)_j\in \{-1,0,1\}$ for $j=1,
2, \ldots, d$, and
\begin{equation}
H_{N, \go} = \gb
  \sum_{n=1}^N \left(1+ \gep \go_n\right) \ind_{\{S_n=0\}}.
\end{equation}
for $\gb \ge 0$ and $\gep \ge 0$. The random variable $\go_1$ is
chosen such that $\bbE[\exp (\gl\go_1)]<\infty$ for every
$\gl\in\R$, and centered. We write $\free (\gb , \gep)$ for $\free$:
 by super--additive arguments $\free$ exists and it is self--averaging (this
 observation is valid for all the models we consider and will not be repeated).
We note that for $\gep =0$ the model can be solved, see e.g.
\cite{cf:G}, and in particular $\free (\gb, 0)=0$ if and only if
$\gb \le \gb_c (d):= -\log (1-\bP(S$ never comes back to $0))$.
Adding the disorder makes this model much more complex: the annealed
bound yields $\free (\gb ,\gep) =0$ if $ \gb \le \gb_c(d)- \log \bbE
\left[\exp (\gep \go_1)\right] =:\tilde {\gb_c}$. It is an open
question whether $\tilde {\gb_c}$ coincides with the quenched
critical value  or not, that is whether $\free (\gb , \gep) =0$
implies $ \gb \le \tilde {\gb_c}$ or not. For references about this
issue we refer to~\cite{cf:AS} and~\cite{cf:Pet}, see however also
the next paragraph:   the model we are considering  can in fact be
mapped to the wetting problem (\cite{cf:AS,cf:G}). Proposition
\ref{th:main} applies to this context with $D_n= \{0\}^\complement$
for every $n$ \cite[Ch. 3]{cf:Feller} and says that one cannot
answer this question via constrained annealed bounds.

\subsubsection{Wetting models in $1+d$ dimensions}
Let $S$ and $\go$ be as in the previous example and
\begin{equation}
\label{eq:Wetting} H_{N, \go} =
\begin{cases}
 \gb \sum_{n=1}^N \left( 1 + \gep \go_n\right) \ind_{\{(S_n)_d=0\}}
 &\text{ if } (S_n)_d\ge 0 \text{ for } n=1,2, \ldots , N \\
 -\infty &\text{ otherwise.}
 \end{cases}
\end{equation}
with $\gb \ge 0$ and $\gep\ge 0$. If one takes  the directed walk
viewpoint, that is if one considers the walk $\{(n,S_n)\}_n$, then
this is a model of a walk constrained above the (hyper--)plane
$x_d=0$ and rewarded $\gb$, on the average, when touching this
plane. If $d=1$ then this is an effective model for a
(1+1)--dimensional interface above a  wall which mostly attracts it.
As a matter of fact in this case there is essentially no loss of
generality in considering $d=1$, since localization is measured in
terms of orthogonal displacements of the walk with respect to
 the wall and we may restrict ourselves  to  this coordinate.
Once again if $\gep =0$ the model can be solved in detail, see e.g.
\cite{cf:G}. Computing the critical $\gb$ and deciding whether the
annealed bound is sharp, at least for small $\gep$, is an unresolved
and disputed question in the physical literature, see e.g.
\cite{cf:FLNO,cf:DHV,cf:TC}. Proposition \ref{th:main} applies with
the choice $D_n =\bbZ ^{d-1}\times \N$.

\subsubsection{Copolymer with adsorption models}
\label{sec:cop} For definiteness choose $S$ to be a one dimensional
simple random walk and take the directed walk viewpoint. Imagine
that the space above the horizontal axis  is filled with a solvent
$A$, while below  there is a solvent $B$. We choose $\go_1 \in
\{A,B\}$ and for example
\begin{equation}
H^{AB}_{N, \go}(S)\, =\, \sum_{n=1}^N \left(a
\ind_{\{\sign(S_n)=+1,\,  \go_n =A\}} + b\ind_{\{\sign (S_n) =-1, \,
\go_n =B\}}+ c\ind_{\{S_n=0\}} \right)
\end{equation}
with $a$, $b$ and $c$ real parameters and $\sign(S_n)=
\sign(S_{n-1})$ if $S_n=0$ (this is just a trick to reward the bonds
rather than the sites). In order to apply Proposition \ref{th:main}
one has to subtract a disorder dependent term, cf.
Remark~\ref{th:rem}: if $a \ge b$ we change the Hamiltonian
\begin{equation}
\label{eq:H1} H_{N, \go}(S)\, :=\, H^{AB}_{N, \go}(S)-
\sum_{n=1}^N a \ind_{\{ \go_n =A\}}.
\end{equation}
without changing the measure $\bP_{N, \go}$  while the free energy
has the trivial shift
 from $\free$ to $\free - a \bbP \left( \go_1=A\right)$.
One can therefore choose $D_n =\bbZ ^{d-1}\times \N$ and Proposition
\ref{th:main} applies. This model has been considered for example in
\cite{cf:ORW}.

Note that if  $c=0$ the model can be cast in a form that has been
considered by a variety of authors (see e.g.
\cite{cf:GHLO,cf:Sinai,cf:AZ,cf:BdH,cf:SSE,cf:TM,cf:Monthus,cf:BG2}):
\begin{equation}
\label{eq:H2} H_{N,\go} (S)\, = \, \gl  \sum_{n=1}^N  \left(\go_n +h
\right) \sign (S_n),
\end{equation}
with $\go$ taking values in $\R$. Once again the Hamiltonian has to
be corrected by subtracting the term $\gl \sum_n (\go_n +h)$ in
order to apply Proposition \ref{th:main}. One readily sees that
\eqref{eq:H1} and  \eqref{eq:H2} are the same model when in the
second case $\go$ takes only the values $\pm 1$, $A=+1$ and $B=-1$,
and $h= (a-b)/(a+b)$, $\gl= (a+b)/4$.

Proposition \ref{th:main} acquires some interest in this context
given the fact that the physical literature is rather split on the
precise value of the critical curve and on whether the annealed
bound is sharp or not, see \cite{cf:BG2} for details on this issue.
In
 \cite{cf:CGG} we present numerical evidence on the
 fact that the annealed curve does not coincide with
 the quenched one, and in view of Proposition~\ref{th:main} this would mean that
 constrained annealing via local functions cannot capture
 the phase diagram of the quenched system.

\subsubsection{Further linear chain models and  observations}
In spite of substantial numerical evidence that in several
instances $\free =0$ but $\freea>0$, we are unaware of an {\sl
interesting} model for which this situation is rigorously known to
happen. Consider however the case $\bbP (\go_1 =+ 1)=\bbP (\go_1 =-
1)=1/2$ and
\begin{equation}
H_{N, \go} (S)=
 \gb \sum_{n=1}^N \left( 1 + \gep \go_n\right) \ind_{\{S_n=n\}},
 \end{equation}
with $\gb$ and $\gep$ real numbers and $S$ the standard simple
symmetric random walk on $\Z$. We observe that Proposition
\ref{th:main} applies to this case with $D_n =\{ n\}^\complement$
and that the model is solvable in detail. In particular $\free (\gb,
\gep) = (\gb - \log 2)\vee 0$, regardless of the value of $\gep$.
The annealed computation instead yields $\freea (\gb , \gep)= (\gb
+\log \cosh(\gep)-\log 2)\vee 0$. Notice in particular that the
critical values of $\gb$, respectively $\log 2$ and $\log 2 -\log
\cosh(\gep)$, differ as long as there is disorder in the system
($\gep \neq 0$).
It is interesting to see in this toy model how the {\sl optimal
choice} of~$A_N$, mentioned at the end of \S~\ref{sec:subintro}, is
rather far from being the empirical average of a local function,
when~$N$ is large.

\smallskip

\begin{rem}\rm
We point out that we restricted our examples only to cases in which
$S$ is a simple random walk, but
 in principle our approach goes through for much more general models,
 like walks with correlated  increments or self--interacting walks, see
 \cite{cf:OTW}  for an example.
 And of course $S_n$ takes  values in $\Z^d$ only for ease of exposition
 and can be easily generalized.
 Another important class of models to which our arguments
 apply is the disordered Poland--Scheraga one \cite{cf:GM2}.
 \end{rem}

 \subsection{The set--up (II): interface pinning models }
 It is natural to wonder whether one can go beyond
 the linear chain set--up. The answer is positive  and
 we give the example of  $(d+1)$--dimensional effective interface models, $d>1$,
 natural generalization of the $(1+1)$--dimensional
 interfaces  considered
 in the previous section.  By this we mean for example
 the case of $S:=\{ S_n \}_{n \in \Z^d}$ with~$S_n \in \R$ and the law of $S$
 is $\bP=\bP_N$:
 \begin{equation}
 \label{eq:d+1free}
 \bP \left( \dd \gp\right)\, \propto \,
 \exp\left(-\frac 12 \sum_{n, n^\prime : \vert n -n^\prime \vert=1} U \left( \gp_n -\gp_{n^\prime}\right)\right)
 \prod_{n \in V_N} \dd \gp_n \prod_{n \in V_N ^\complement} \gd _0 (\dd \gp _n)
 ,
 \end{equation}
 where $V_N= [-N/2, N/2]^d \cap \Z^d$ and
  $U(\cdot)$ is a measurable function such $\lim_{r \to \pm \infty} U(r)=+\infty$
 sufficiently rapidly to make the right--hand side of \eqref{eq:d+1free} integrable
 (note that we may assume $U(\cdot)$ to be even). As a matter of fact, in order to
 have a treatable model one has to restrict rather strongly the choice of $U(\cdot)$:
 interface models are extremely challenging even without introducing pinning
 potentials (or, of course, disorder). Connected to that  is also the reason why
 we have chosen the continuous set--up for interface models: discrete models
 are even more challenging \cite{cf:Grev}.

 The disorder in the system this time is given by an IID {\sl field} $\go := \{ \go_n \}_{n \in \Z ^d}$
 and $H_{N,\go}(S)$ depends only upon $S_n$ with $n \in V_N$: $\go_0$ takes once again
 values in~$\Gamma$.  The definition  \eqref{eq:1}
 of $\bP_{N, \go}$ is unchanged and the Basic Hypothesis varies in the obvious way, that is
 we assume that there exists $\{ D_n\}_{n \in \Z^d}$
 such that
 \begin{equation}
 \label{eq:BHd}
 \lim_{N \to \infty} \frac 1{N^d}
 \log \bP \left( S_n \in D_n \text{ for } n \in V_N \right)\, =\, 0,
 \end{equation}
and such that $H_{N, \go } (S)=0$ if $S_n \in D_n$ for every  $n \in
V_N$. Like for linear chains we assume the existence of the quenched
free energy, that is of the $L^1(\bbP(\dd \go))$ and $\bbP( \dd
\go)$--a.s.   limit of the sequence $\left\{ N^{-d} \log Z_{N, \go}
\right\}_N$ and like in the linear chain case  we have $0 \le f \le
\tilde f$, where $\tilde f$ is again the annealed free energy
defined in analogy with \eqref{eq:J}.

\smallskip

 The punch--line of this section is that Proposition~\ref{th:main}
 holds in this new set--up and it is proven exactly in the same way:

 \medskip

\begin{proposition}
\label{th:main-d} If $\freea >0$ then for every local bounded
measurable function $F: \Gamma ^{\Z^d}
 \longrightarrow \R$ such that $\bbE \left[ F (\go) \right]=0$
one has
\begin{equation}
\label{eq:main-d} \liminf_{N \to \infty} \, \frac{1}{N^d} \, \log
\bbE \bE \Bigg[ \exp \Bigg( H_{N, \go}(S) + \sum_{n\in \gL_N}
F(\theta_n \go)\Bigg)\Bigg]\, > \, 0.
\end{equation}
\end{proposition}

\medskip

 In order to give  examples of applications we may consider the
 $d+1$ dimensional  model of random rewards and
 penalties near the origin, that is the case of
 \begin{equation}
\label{eq:RRP} H_{N, \go} \, =\,  \gb \sum_{n\in V_N} \left( 1 +
\gep \go_n\right) \ind_{\{S_n\in(-1,1)\}},
\end{equation}
but one can write natural straightforward generalizations of the
wetting models and of the copolymer with adsorption. The Basic
Hypothesis in all these cases is a probability estimate on what is
known as an {\sl entropic repulsion event}, that is, for example,
the event that $S_n \ge 1$ for every $n \in V_N$ and one can for
example show that such a probability is bounded below by
 $\exp\left( -c N^{d-1}\right)$, $c>0$,
if $U(\cdot)$ is $C^2$ and $ \inf_r U^{\prime \prime} (r)  >0$, see
 \cite{cf:Grev} and references therein.
So in this case one may apply Proposition~\ref{th:main} to conclude
that one cannot improve on the annealed bound by constraining via
local functions.

Two comments, of opposite spirit, are however in order (for details
see the lecture notes \cite{cf:Grev}):
\smallskip
\begin{enumerate}
\item
The Basic Hypothesis requires a substantially weaker estimate and it
is reasonable to expect that one is
 able to  verify it in greater generality.
\item\rule{0pt}{14pt}The understanding of
 the associated deterministic models ($\gep=0$ for random rewards and wetting models
 and the annealed models in general)
 is still extremely  partial. Somewhat satisfactory results are available for
 quadratic $U(\cdot)$, that is $\bP$ is Gaussian, but even in this case
 one has  to give up the precise estimates
 available for the linear chain case (like computing exactly $\gb_c$)
 and basic questions are still open. So the application
 of  Proposition~\ref{th:main-d}, while being relevant on a conceptual level,
  yields a result that has little quantitative content.
 \end{enumerate}

\bigskip

\section{On zero free energy and  null potentials}
\label{sec:excocycles}

In this Section $d\ge 1$. Let $\{\go_n\}_{n\in\Z^d}$ be an IID
family of random variables under the probability measure $\bbP$,
taking values in $\Gamma=\R$. The law of $\go_1$  is denoted by
$\nu$.

We are interested in the family $A = \{A_N\}_{N\in\N}$
 of empirical averages of a  local function $F$, that is
\begin{equation} \label{eq:hamiltonian}
A_N (\go)\, = \, \sum_{n \in V_N} F\left( \theta_n {\go}\right)\,,
\end{equation}
where $F: \gG ^{\Z ^d} \to \R$ depends only on the variables indexed
by a finite set $\gL \subset \Z^d$, that is $F(\go)=F(\go^\prime)$
if $\go_n =\go ^\prime_n$ for every $n \in \gL$. Notice that, by
standard (super--additivity) arguments, the limit
\begin{equation} \label{eq:L}
L(F) \, := \, \lim_{N\to \infty} \frac 1{N^d} \log
\bbE\left[\exp\left(  A_N(\go)\right) \right],
\end{equation}
exists. Moreover, by Jensen's inequality, $ L(F) \ge \bbE \left[ F (
\go)\right]$.

\smallskip

We will prove the following:

\bigskip
\begin{proposition}
\label{th:main2} Assume that $\bbE\left[ F( \go )\right]=0$. If
$L(F)=0$, then
\begin{equation} \label{eq:main2}
    \lim_{N \to \infty} \, \frac{1}{N^d} \, \sup_\go \vert A_N (\go)\vert =0 \,.
\end{equation}
\end{proposition}
\bigskip

Of course, since the result is uniform in $\go$, the proposition
covers also the linear chain set--up, where one considers
$\theta_{\lfloor N/2\rfloor +1}V_N$ rather than $V_N$.

\medskip
\noindent {\it Proof.} We consider the {\sl potential}, in the sense
of \cite[Def. (2.2)]{cf:Georgii}, $\Phi  := \left\{
\Phi_B\right\}_{B \subset \Z^d}$ defined by
\begin{equation}
\Phi_B (\go) \, = \,
\begin{cases}
F\left( \theta_{-n} \go\right) &\text{ if there exists } n \text{ such that } \theta_n B= \gL ,\\
0 & \text{ otherwise}.
\end{cases}
\end{equation}
Let $\nu $ be the {\sl single spin} reference measure \cite[Def.
(2.9)]{cf:Georgii} and let us set
\begin{equation}
Z_{N}^\Phi (\go) \, :=\, \int \exp\left(H _N^\Phi (\gs)\right) \prod
_{n\in V_N} \nu (\dd \gs_n) \prod _{n\in V_N^\complement}
\delta_{\go_n} (\dd \gs_n),
\end{equation}
with $H _N^\Phi (\gs):=\sum_{B: B\cap V_N \neq \emptyset}\Phi_B
(\gs)$. Note that $A_N(\cdot)$ differs from $H _N^\Phi (\cdot)$ only
by boundary terms so that~$\sup_{\go} |A_N(\go) - H_N^\Phi(\go)| \le
C N^{d-1}$ for some~$C>0$ (we recall that $F(\cdot)$ is bounded).
Therefore it suffices to show that~\eqref{eq:main2} holds with
$A_N(\cdot)$ replaced by~$H _N^\Phi (\cdot)$.

\smallskip

Let us consider  the $\theta$--invariant Gibbs measure $\mu$
associated to the potential~$\Phi$, the existence of which is
established in a standard way by taking infinite volume limits with
periodic boundary conditions (if $\nu$ has unbounded support
tightness follows from the fact that $F(\cdot)$ is bounded). By
\cite[Theorem (15.30)]{cf:Georgii} the {\sl relative entropy
density} of $\nu ^{ \infty}$ ($\nu^\infty (\dd \go):= \prod_{n \in
\Z^d} \nu (\dd \go_n)$)
 with respect to
$\mu$ exists and can be written as
\begin{equation}
\lim_{N\to \infty} \, \frac 1{N^d} \, \cH_{V_N}\Big(  \nu ^{\infty}
\,\big|\, \mu \Big)
 \, =\, \lim_{N\to \infty} \, \frac 1{N^d}  \log
Z_{N}^\Phi (\go)  \, - \, \int F (\go) \, \nu^{\infty} (\dd \go),
\label{eq:relentest}
\end{equation}
where $ \cH_{V_N}(\nu ^{\infty} | \mu)$ is the relative entropy of
$\nu^\infty$ with respect to $\mu$, when both measures are
restricted  to the $\gs$--algebra generated by the variables
$\{\go_n\}_{n \in V_N}$. We have of course used  the standard
definition of relative entropy, $\cH(\mu_1\vert \mu_2)= \int \log
(\dd \mu_1/\dd \mu_2) \dd \mu_1$ for $\mu_1$ and $\mu_2$ two
probability measures with $ \mu_1$ absolutely continuous with
respect to $ \mu_2$. A last remark on formula \eqref{eq:relentest}
is that it holds for any choice of $\go$: this is just the
independence of the free energy on boundary conditions. This
independence may be seen directly since $\log
(Z_N^{\Phi}(\go)/Z_N^{\Phi}(\go^\prime))= O(N^{d-1})$ uniformly in
$\go$ and $\go^\prime$ and this implies also that
 the first
term in the right--hand side of \eqref{eq:relentest} may be replaced
by $L(F)$.

Notice now that both terms
 in the right--hand side of \eqref{eq:relentest}
are zero, respectively by the hypotheses $L(F)=0$ and
$\bbE[F(\go)]=0$, and therefore, as a consequence of
 the Gibbs variational principle \cite[Theorem (15.37)]{cf:Georgii}, $\nu^{\infty}$ is a Gibbs measure
with the same specification of $\mu$, but of course $\nu^{\infty}$
is the Gibbs measure with potential $\Phi^{(0)}$ identically equal
to zero and single spin measure $\nu$. This means that $\Phi
-\Phi^{(0)}(=\Phi)$ is a {\sl negligible} potential, that is
\cite[Theorem (2.34)]{cf:Georgii} the function
\begin{equation}
 \sum_{B: B\cap V_N \neq \emptyset}
\left( \Phi_B (\go) - \Phi^{(0)}_B (\go) \right)
 \end{equation}
does not depend on the variables $\go_n$ for~$n \in V_N$. We can
write
\begin{equation}
\begin{split}
\cH_N^\Phi(\go)\, =\, \sum_{B: B\cap V_N \neq \emptyset} \Phi_B
(\go) \, &=\,  \sum_{B: B\subset V_N } \Phi_B (\go)+
 \sum_{B: B\cap V_N \neq \emptyset, \, B\not\subset V_N  } \Phi_B (\go)\\
 &=:\, I_N(\go) + R_N(\go),
 \end{split}
\end{equation}
and since $\cH^\Phi_N(\go)$ does not depend on the $\go_n$'s for
$n\in V_N$ we may change  in the right--hand side the  configuration
$\go$ with $\tilde\go $ defined by setting
 $\tilde \go_n = \go_n$ for $n \in V_N^\complement$
and $\go_n=c$, $c$ an arbitrary fixed constant, for $n \in V_N$.
Therefore, in random variable terms, we have
\begin{equation}
\cH_N^\Phi(\go)\, =\, c_N+ R_N(\tilde{ \go} ),
\end{equation}
with $c_N= I_N(\tilde {\go})$ (notice that it is not random and it
depends only on the choice of $c$). From the immediate estimate
$\sup_\go \vert R_N (\go) \vert \le C N^{d-1}$ for some~$C=C(F)>0$
it follows that for all~$\go$
\begin{equation}
    c_N -  C N^{d-1} \,\le\, \cH_N^\Phi(\go)\, \le \, c_N +  C N^{d-1},
\end{equation}
and the hypothesis $L(F)=0$ yields immediately $\lim_{N\to \infty}
c_N/N^d=0$. Therefore
\begin{equation}
 \sup_\go \left \vert \cH_N^\Phi(\go)\right\vert \, \le \, c_N + C N ^{d-1} \,= \, o(N^d),
\end{equation}
and the proof is complete.

\section*{Acknowledgments}
The contributions of the referees to this note have been very
important. The first version of this paper was restricted to linear
chain models and the proof was based on cocycles and
Perron--Frobenius theory (this version is still available in F.C.'s
Ph.D. Thesis). We owe the approach in Section 2 of the present
version to the suggestion of one of the referees. We would also like
to thank T. Bodineau, E. Orlandini and F. L. Toninelli for helpful
discussions.



\begin{thebibliography}{99}

\bibitem{cf:AZ}
S. Albeverio and X. Y. Zhou, \textit{Free energy and some sample
path properties of a random walk with random potential}, J. Statist.
Phys. {\bf 83 } (1996),  573--622.

\bibitem{cf:AS}
 K. S. Alexander and V. Sidoravicius,
 \textit{Pinning of polymers and interfaces by random potentials},
 preprint (2005). Available on: arXiv.org e-Print archive: math.PR/0501028



\bibitem{cf:BG2}
T. Bodineau and G. Giacomin, \textit{On the localization transition
of random copolymers near selective interfaces}, J. Statist. Phys.
{\bf 117} (2004), 801--818.

\bibitem{cf:BdH}
E. Bolthausen and F. den Hollander, \textit{Localization transition
for a polymer near an interface}, Ann. Probab.  {\bf 25}  (1997),
1334--1366.

\bibitem{cf:CGG} F. Caravenna, G. Giacomin and M. Gubinelli,
\textit{A numerical approach to copolymers at selective interfaces},
preprint (2005), available on hal.ccsd.cnrs.fr


\bibitem{cf:DHV}
B. Derrida, V. Hakim and J. Vannimenus, \textit{Effect of disorder
on two--dimensional wetting}, J. Statist. Phys. {\bf 66} (1992),
1189--1213.





\bibitem{cf:add1}
A. C. D. van Enter, C. K\"ulske and C. Maes, \textit{Comment on:
Critical behavior of the randomly spin diluted 2D Ising model: A
grand ensemble approach (by R. K\"uhn)}, Phys. Rev. Lett. {\bf 84}
(2000), 6134.

\bibitem{cf:Feller}
W.~Feller, \textit{An introduction to probability theory and its
applications}, Vol. I, Third edition, John Wiley \& Sons, Inc., New
York--London--Sydney, 1968.


\bibitem{cf:FLNO}
G. Forgacs, J. M. Luck, Th. M. Nieuwenhuizen and H. Orland,
\textit{Wetting of a disordered substrate: exact critical behavior
in two dimensions}, Phys. Rev. Lett. {\bf 57} (1986), 2184--2187.

\bibitem{cf:GM2}
T. Garel and C. Monthus, \textit{Numerical study of the disordered
Poland-Scheraga model of DNA denaturation},
 J. Stat. Mech., Theory and Experiments (2005), P06004.

\bibitem{cf:Georgii}
H.--O.~Georgii,  \textit{Gibbs measures and phase transitions}, de
Gruyter Studies in Mathematics, { \bf 9}. Walter de Gruyter \& Co.,
Berlin, 1988.

\bibitem{cf:G}
G.~Giacomin, \textit{Localization phenomena in random polymer
models}, preprint (2004), available on
www.proba.jussieu.fr/pageperso/giacomin/pub/publicat.html

\bibitem{cf:Grev}
G.~Giacomin, \textit{Aspects of statistical mechanics of random
surfaces}, unpublished manuscript, notes of Lectures given at the
IHP, Paris, in the fall 2001, available on
www.proba.jussieu.fr/pageperso/giacomin/pub/publicat.html

\bibitem{cf:GT}
G.~Giacomin and F.~L.~Toninelli, \textit{Estimates on path
delocalization for copolymers at interfaces},
 Probab. Theory Rel. Fields (online first, May 2005).

\bibitem{cf:GHLO}
T. Garel, D. A. Huse, S. Leibler and H. Orland, \textit{Localization
transition of random chains at interfaces}, Europhys. Lett. {\bf 8}
(1989), 9--13.


\bibitem{cf:Kuhn}
R. K\"uhn, \textit{Equilibrium ensemble approach to disordered
systems I: general theory, exact results}, Z.~Phys.~B (1996),
231--242.

\bibitem{cf:add2}
C. K\"ulske, \textit{Weakly Gibbsian Representations for joint
measures of quenched lattice spin models}, Probab. Theory Rel.
Fields {\bf 119} (2001), 1--30.



\bibitem{cf:LDMF99}
P. Le Doussal, C.  Monthus and  D. S.  Fisher, \textit{Random
walkers in one-dimensional random environments: exact
renormalization group analysis}, Phys. Rev. E (3)  {\bf 59}  (1999),
4795--4840.


\bibitem{cf:Monthus}
C. Monthus, \textit{On the localization of random heteropolymers at
the interface between two selective solvents}, Eur. Phys. J. B {\bf
13} (2000), 111--130.

\bibitem{cf:Morita}
T. Morita, \textit{Statistical mechanics of quenched solid solutions
with application to magnetically dilute alloys}, J. Math. Phys. {\bf
5} (1966), 1401--1405.

\bibitem{cf:ORW}
E. Orlandini, A. Rechnitzer and S. G. Whittington, \textit{Random
copolymers and the Morita approximation: polymer adsorption and
polymer localization}, J. Phys. A: Math. Gen.  {\bf 35} (2002),
7729--7751.

\bibitem{cf:OTW}
E. Orlandini, M. C. Tesi and S. G. Whittington, \textit{A
self--avoiding model of random copolymer adsorption}, J. Phys. A:
Math. Gen.  {\bf 32} (1999), 469--477.

\bibitem{cf:Pet}
N. Petrelis, \textit{Polymer pinning at an interface}, preprint
(2005). Available on: arXiv.org e-Print archive: math.PR/0504464

\bibitem{cf:Sinai}
Ya. G. Sinai, \textit{A random walk with a random potential}, Theory
Probab. Appl. {\bf 38} (1993),  382--385.

\bibitem{cf:SSE}
S. Stepanow,  J.-U. Sommer and   I. Ya. Erukhimovich,
\textit{Localization transition of random copolymers at interfaces},
Phys. Rev. Lett. {\bf 81} (1998), 4412--4416.

\bibitem{cf:TC}
L.--H. Tang and H. Chat\'e, \textit{Rare--event induced binding
transition of heteropolymers}, Phys. Rev. Lett. {\bf 86} (2001),
830--833.



\bibitem{cf:TM}
A. Trovato  and A. Maritan, \textit{A variational approach to the
localization transition of heteropolymers at interfaces}, Europhys.
Lett. {\bf 46} (1999), 301--306.

\end{thebibliography}
\end{document}